\documentclass[11pt]{amsart}
\hoffset= - 0.75 in

\usepackage{amsfonts, amsmath, amssymb, amsthm}
\usepackage{latexsym, graphicx, color, graphics}
\newtheorem{thm}{Theorem}[section]

\newtheorem{defn}[thm]{Definition}

\newtheorem{cor}[thm]{Corollary}
\newtheorem{prop}[thm]{Proposition}

\newcommand{\RR}{\mathbb{R}}

\newcommand{\TP}{\mathbb{TP}}

\begin{document}

\title{Tropical secant varieties of linear spaces}
\author{Mike Develin}
\address{Mike Develin, American Institute of Mathematics, 360 Portage 
Ave., Palo Alto, CA 94306-2244, USA}
\date{\today}
\email{develin@post.harvard.edu}

\begin{abstract}
In this paper, we investigate tropical secant varieties of ordinary linear 
spaces. These correspond to the log-limit sets of ordinary toric 
varieties; we show that their interesting parts are combinatorially isomorphic to a certain 
natural subcomplex of the complex of regular subdivisions of a corresponding point 
set, and we display the range of behavior of this object. We also use this characterization to 
reformulate the question of determining Barvinok rank into a question 
regarding regular subdivisions of products of simplices.
\end{abstract}

\maketitle

\section{Introduction}
The tropical semiring is given by the real numbers $\RR$, together with 
the operations of tropical addition $\oplus$ given by $a\oplus b = 
\text{min}\,(a,b)$ and tropical multiplication $\otimes$ given by 
$a\otimes b = a+b$. As in ordinary geometry, we consider the tropical 
semimodule $\RR^{d+1}$, as well as the corresponding tropical projective 
space $\TP^{d} = \RR^{d+1}/(1,\ldots,1)$ given by modding out by tropical 
scalar multiplication. There has been a recent spate of work in tropical 
geometry, as well as the use of tropical geometry to solve problems 
arising in ordinary geometry. In particular, theories of tropical 
convexity~\cite{DS}, tropical polytopes~\cite{Joswig}, tropical linear 
spaces~\cite{SS}, tropical linear algebra~\cite{DSS}, tropical 
geometry~\cite{RGST}, and tropical algebraic geometry~\cite{Mikh} have all 
been burgeoning.

In this paper, we consider the {\em k-th tropical secant variety} of an
ordinary linear subspace $L$ in tropical projective space $\TP^d$, defined 
by 
\[
S^k(L) := \{v_1\oplus v_2\oplus\cdots\oplus v_{k+1}\,\mid\, x_i\in L\}.
\]
These tropical secant varieties correspond to the log-limits of ordinary 
toric varieties. To be precise, the image of an ordinary toric variety in 
the variables $x_1,\ldots,x_n$ under the logarithm map is a linear space 
in the variables $\text{log}\,x_1, \ldots, \text{log}\,x_n$; the dimension of the linear space is the 
same as the dimension of the toric variety. For $x_i$ and 
$y_i$ large, we have $\text{log}\,(x_i+y_i)\sim 
\text{max}\,(\text{log}\,x_i + {\rm log}\,y_i)$, and so the log-limit of 
the ordinary secant variety corresponds to the tropical secant variety 
(the max-plus and min-plus semirings are isomorphic.)

Another application of tropical secant varieties of linear subspaces is to
the Barvinok subcomplex. A matrix of size $d\times n$ has {\em Barvinok
rank} at most $k$ if it is expressible as the tropical sum of $k$
tropically rank-one matrices, where a tropically rank-one matrix $M$ is
one satisfying $m_{ij}+m_{kl} = m_{il}+m_{jk}$ for all
$i,j,k,l$~\cite{DSS}. These tropically rank-one matrices form an ordinary
linear subspace, and the matrices of Barvinok rank at most $k$ comprise
the $k$-th secant variety of this linear subspace. A consequence of our 
work is that the interesting component of this space is a subcomplex of 
the secondary polytope of the product of simplices $\Delta_{d-1}\times 
\Delta_{n-1}$; we also use our results to provide an intuitive algorithm 
for determining Barvinok rank.

Our aim is to develop a general theory of these tropical secant varieties
of ordinary linear spaces. In Section~\ref{maintheorem}, we present and
prove our main theorem, which states that the interesting parts of tropical secant varieties of a 
linear space, which we 
call {\em tropical secant complexes}, 
are certain natural subcomplexes of the complex of regular subdivisions
of a corresponding point configuration. In Section~\ref{1-dim}, we use this
representation to prove that the $k$-th secant variety of any generic line
in $d$-space is equal to the cone from a line over the complex of lower
faces of the cyclic polytope $C(2k,d-1)$. In Section~\ref{examples}, we
compute the first secant variety of a diverse set of two-dimensional
examples, including an example where the corresponding complex is not pure and one where it is not 
contractible.
In Section~\ref{barvinok}, we apply our theory to the
case of $d\times n$ matrices of Barvinok rank $k$, showing that these
complexes are certain subcomplexes of the secondary polytope of
$\Delta_{d-1}\times \Delta_{n-1}$. A corollary of this is that the complex
of matrices of Barvinok rank two is pure.

\section{Tropical secant complexes}~\label{maintheorem}

In this section, we prove the following main theorem allowing us to express tropical secant varieties 
of linear spaces as subcomplexes of the fan of regular subdivisions of an associated point 
configuration.

\begin{thm}\label{thm-main}
Let $L\subset \TP^{n-1}$ be the (ordinary) linear subspace of dimension $d$ generated by the $d$ rows 
of the associated matrix $M_L$. Let $V_L=\{v_1,\ldots,v_n\}$ be the $n$-point configuration in $\RR^d$ 
given by the 
columns of $M_L$. Then a vector $x=(x_1,\ldots,x_n)$ is in the $k$-th tropical secant variety of $L$ if 
and only 
if the upper envelope of the polytope formed by the height 
vector $x$ has $k+1$ facets whose union contains each point of $V_L$.
\end{thm}

\begin{proof}
Let the rows of $M_L$, the generators of $L$, be denoted by $L_1,\ldots,L_k$, and let the corresponding
coordinates of $V_L$ be $z_1,\ldots,z_k$. A point $x$ is in the $k$-th secant variety of $L$ if and
only if there exist points $y^1, \ldots, y^{k+1}$ in $L\oplus (1,\ldots,1)\RR$ such that $x_i =
\text{min}_j\:(y^j_i)$. Consider any point $y^j\in L$. This is some linear combination of the rows of
$L$, plus a constant:  $y^j = a + \sum_{k=1}^d c_k L_k$. The corresponding height function on the
points $v_i$ is given simply by $f^j = a+\sum c_k z_k$; in other words, the 
heights given by $y^j$ to 
$V_L$ are given by the value of some affine functional on $\RR^d$ at those points.

Therefore, a point $x$ is in the $k$-th secant variety if there exist $k+1$ affine functionals on
$\RR^d$ whose pointwise minima at the points of $V_L$ give the coordinates of $x$. Now, consider the
regular subdivision of $V_L$ given by the upper envelope of the height vector $x$ as in the statement
of the theorem. If these $k+1$ affine functionals exist, then each one 
defines an upper face of this
polytope, since for each functional $f^j$ and each point $v_i$ we have $x_i \le f^j(v_i)$ (since $x_i$
is the coordinate-wise minimum of the $f^j$. Furthermore, for each $v_i$, $x_i$ is equal to some
$f^j(v_i)$, and so $v_i$ is in the face defined by $f^j$. Therefore, these faces together contain each
point of $V_L$, and thus so do a set of $k+1$ facets containing them.

Conversely, suppose that the height vector $x$ induces a regular subdivision with $k+1$ upper facets
$F_1,\ldots,F_{k+1}$ such that each $v_i$ is contained in one of them. Each facet $F_j$ is contained in
some hyperplane $H_j$ in the lifted $\RR^d$ with height vectors, and as before this hyperplane
corresponds to an affine functional on $\RR^d$. If this functional is $a+\sum c_k z_k$, then define the
point $y^j$ via $y^j = a + \sum_{k=1}^d c_k L_k$. We claim that $y^1\oplus\cdots\oplus y^{k+1} = x$. We
need to check that they agree in each coordinate, which is the same thing as saying that the height
vector $x$ is the coordinate-wise minimum of the height vectors $y^j$. However, since $F_j$ is an upper
facet for each $j$, for each $i$ we have $x_i \le y^j_i$, and since $x_i$ is contained in one of the 
facets, equality is achieved for some $j$. Therefore, $x$ is in the $k$-th secant variety of $L$ as 
desired.
\end{proof}

It is worth noting here that we picked an arbitrary basis for $L$. However, picking a different basis 
yields an affinely isomorphic point configuration, so as must be the case we can pick any basis to fill 
out the matrix $M_L$.

Since each regular subdivision corresponds to a polyhedral cone of height vectors,
Theorem~\ref{thm-main} gives us a decomposition of the tropical secant variety.

\begin{cor}\label{complex}
The $k$-th tropical secant variety of a linear subspace $L$ is a cone from $L$ over a polytopal
complex, which we call the {\em k-th tropical secant complex} of $L$. The faces of this polytopal 
complex correspond to regular subdivisions of $L$ in which there
exist $k+1$ facets containing all of the points, with a face $F$ containing a face $G$ if the regular
subdivision associated to $F$ refines the one associated to $G$.
\end{cor}

\begin{proof}
First of all, consider the case of a regular subdivision not using all of the vertices of the point
configuration.  In order for a height vector to be valid under Theorem~\ref{thm-main}, the lift of any
unused point $x$ must lie on the lifted face whose interior contains $x$, and in particular the height
of $x$ is uniquely specified by the vertices of the point configuration, with the space of height
vectors inducing that regular subdivision being affinely isomorphic to the space of height vectors on
the used vertices inducing the regular subdivision on that point subconfiguration.

This space is clearly a polyhedral cone, with inequalities on the heights given by the set of upper
facets; in other words, the inequalities are given by picking a vertex, picking an upper facet, and 
noting that the point's height must be less than that of the corresponding facet-defining hyperplane. 
Setting one of these inequalities to an equality corresponds to the case where a point off a facet is 
moved onto the facet, which has the effect of coarsening the regular subdivision, so these cones fit 
together as stated in the corollary. The original cone is in the $k$-th secant variety if and only if 
the corresponding regular subdivision has $k+1$ facets whose union is all of the points; this property 
is preserved under coarsening, so every face of any cone in the $k$-th secant variety is also in the 
$k$-th secant variety. Therefore, the $k$-th secant variety is a valid polyhedral complex.

To check that it is a cone from $L$ over a polytopal complex, we need first to check for each cone that
if $x\in L$, then $x$ is a cone point. However, as before, the height vector of $x\in L$ is simply an
affine functional on the space containing the configuration $V_L$. Taking a nontrivial linear
combination of $x$ with any height vector $y$ simply performs an affine transformation on the heights,
which does not change the induced regular subdivision and thus does not change the cone of the height 
vector as desired.

The other step in checking that each cone is a cone from $L$ of a polytope is to check that the only 
lineality in each cone is in fact $L$. Indeed, consider any vector $y\notin L$. Because $y$ is not in 
$L$, it is not an affine function on the configuration $V_L$, so $y$ does not induce the trivial 
subdivision. It immediately follows that $y$ and $-y$ induce different subdivisions, indeed 
subdivisions with no common refinement. But given any $x$ in the cone, for large enough $k$, $x+ky$ 
refines $y$ and $x-ky$ refines $-y$, so it is impossible that these two are the same regular 
subdivision, and hence they cannot both be in the original cone. Therefore, $y$ cannot be in the 
lineality space of the cone as desired.
\end{proof}

In the convex case, this gives us something very akin to the secondary polytope. Given a polytope $P$,
the {\em secondary polytope}~\cite{Z} is a polytope whose face poset is the poset of regular
subdivisions, ordered by refinement, so that the vertices of $P$ are the regular triangulations. If
$V_L$ is in convex position, and $k$ is large, all regular subdivisions satisfy the condition of
Corollary~\ref{complex}. The condition for face inclusion is dual to the inclusion in the secondary
polytope, so the $k$-th secant variety is a cone from $L$ over the dual of the secondary polytope. In 
particular, it will be all of $\RR^n$ in this case. Indeed, this is an if and only if.

\begin{cor}
The $\infty$-th secant variety of a linear subspace $L$ is all of $\RR^n$ if and only if the 
corresponding point configuration $V_L\subset \RR^d$ is in convex position, i.e. if every point of 
$V_L$ is a vertex of $\text{conv}\,(V_L)$. 
\end{cor}

\begin{proof}
A point $x$ is in the $\infty$-th secant variety if and only if its height vector corresponds to a lift 
of $V_L$ 
where the union of all of the upper facets contains all of the lifted points. If $V_L$ is in convex 
position, every height vector will have this property, since every point will be in some upper facet. 
If not, then some point can be written as an affine combination of the other points, $v_i = \sum_{j\neq 
i} c_jv_j$ with $\sum c_j = 1$. Then if $x_i < \sum_{j\neq i} c_j x_j$, the lifted point $v_i$ will not 
be in the upper envelope of the convex hull of all of the lifted points, and so it will be in no upper 
facet. Therefore, any $x$ satisfying this condition will not be in the $\infty$-th secant variety of 
$L$.
\end{proof}

We can easily compute the dimension of the $k$-th tropical secant variety in the case where the linear 
subspace is {\em generic} in the sense that the point configuration $V_L$ is in general position, i.e. 
has no $d+1$ points lying in an affine subspace of dimension $d-1$.

\begin{prop}\label{cxdim}
Suppose that the $d$-dimensional linear subspace $L$ is generic.
Then the $k$-th secant variety of $L$ is a complex of dimension
$\text{min}\,((k+1)(d+1)-1, n)$.
\end{prop}

\begin{proof}
By Corollary~\ref{complex}, the faces of the $k$-th tropical secant variety of $x$ correspond to 
regular subdivisions of the $n+1$-point configuration $V_L$ in $\RR^d$ in which there exist $k+1$ 
facets whose union contains all of the points. The condition of the theorem implies that $V_L$ consists 
of $n+1$ distinct points. Since the tropical secant variety is a subset of $\TP^n$, it obviously must 
have dimension at most $n$. 

First, suppose that $(k+1)(d+1)-1\le n$. Pick a generic linear functional; this orders the points in 
$V_L$
in some order, $v_1,\ldots,v_{n+1}$.  Then by examining the values of this linear functional at these
points, it is evident that the $d$-polytopes formed by the convex hulls of the sets
\[
\{v_1,\ldots,v_{d+1}\}, \{v_{d+2}, \ldots, v_{2d+2}\}, \ldots, \{v_{1+(k-1)d}, \ldots, 
v_{1+(k-1)(d_1)-1}\}, \{v_{1+(k-1)(d+1)}, \ldots, v_n\}
\]
do not intersect. These are $k$ $d$-simplices and one facet which is not a simplex (they have 
dimension $d$ since the point configuration $V_L$ is in general position.) Since these facets 
do not intersect, there exists a regular subdivision containing all of them. Take a regular subdivision 
from this nonempty set which is as fine as possible. We claim that the corresponding cell of the 
$k$-th tropical secant variety has dimension $(k+1)(d+1)-1$. Indeed, consider a height vector in its 
relative 
interior. For each of the $k+1$ $d$-polytopes in our set, we have $d+1$ degrees of freedom for the 
heights corresponding to those vertices, which gives us $(k+1)(d+1)$ degrees of freedom which we can move in without 
changing the subdivision. Furthermore, in the relative interior of this cell, wiggling these points can 
only refine the rest of the induced subdivision, and since the cell was chosen to be maximal with 
respect to refinement, it also leaves the remainder of the subdivision fixed. Therefore, the cell of 
height vectors has  dimension equal to $(k+1)(d+1)$, which upon projectivization yields a cell of 
dimension $(k+1)(d+1)-1$ in the $k$-th tropical secant variety.

On the other hand, any cell of the $k$-th tropical secant variety corresponds to a regular 
subdivision with $k+1$ 
facets whose union is all of $V_L$. Each of these facets is $d$-dimensional, and so in order for a 
height vector to be in this cell, the coordinates of the height vector which correspond to each facet 
must lie in a $d+1$-dimensional space. Since each coordinate corresponds to at least one facet, we have 
at most $(k+1)(d+1)$ degrees of freedom in the cell, and so the maximal dimension of any cell in the 
tropical secant variety is in fact $(k+1)(d+1)-1$ as desired.

Finally, suppose $(k+1)(d+1)-1\ge n$. We need to prove that there exists some cell of dimension $n$ in 
the $k$-th tropical secant variety of $L$. As before, we find a linear 
functional, which 
orders the points $v_1,\ldots,v_{n+1}$. For our (at most) $k+1$ facets, we take as many $d+1$'s as 
possible, until we are left with some remainder less than $d+1$; for the final facet, we take a simplex 
including this leftover subset and some points from the last set of $d+1$, and as in the first case we 
find a triangulation including these facets and as fine as possible given that constraint. Then by the 
same logic as before, it immediately follows that we have $n$ degrees of freedom which we can move in 
without changing the subdivision in question.
\end{proof}

This maximal dimension is what one would expect; given some generic $d$-dimensional shape in ordinary
space, the $k$-th secant variety consists of the union of all points in $k$-planes spanned by $k+1$
points from the shape. The number of degrees of freedom of this space should be as follows: for each of
the $k+1$ points, we have $d$ degrees of freedom, and we have a final $k$ degrees of freedom for
picking the point once the plane is fixed, for a total of $(k+1)d+k = (k+1)(d+1)-1$, just as in the 
tropical secant case. The first tropical secant variety is especially nice.

\begin{prop}\label{first-pure}
If the $d$-dimensional linear space $L$ is generic, then the first tropical secant variety of $L$ is 
a pure complex.
\end{prop}

\begin{proof}
By Proposition~\ref{cxdim}, the dimension of the first tropical secant variety of $L$ is ${\rm 
min}\,(2d+1, n)$. We 
need to show that every regular subdivision which contains two facets whose union contains each point 
in $V_L$ can be refined to one with this many degrees of freedom. First, suppose $2d+1<n$. If our 
two facets do not intersect, then as in the proof of the previous proposition, we take a regular 
subdivision containing the provided one which is as fine as possible. The corresponding cell will then 
have $d+1$ degrees of freedom for each facet for a total of $2d+1$ (after projectivization.)

If the two facets do not intersect, by the previous case, it suffices to show that this complex can be
refined to one with two non-intersecting spanning facets. Suppose they intersect in an $r$-face;  we
will induct on $r$. Since $V_L$ is in general position, this face must be a simplex with $r+1$ points.
Since $n>2d+1$, and we have $n+r+1$ points with multiplicity, one of the two facets must have at least
$d+2$ points. Take one overlap point out of it; the remaining points form a full-dimensional convex
hull.  Then we can refine the previous subdivision by breaking this facet; rigorously, what we are 
doing is lowering the height of the shared vertex by an infinitesimal amount and adjusting the heights 
of facet $F_1$ so that the lifted points still share a hyperplane. This process has the effect of 
taking the shared vertex out of the upper facet-defining hyperplane of $F_2$; the new, refined 
subdivision (which may be strictly finer in other places as well) refines the original one and has two 
facets whose union is all of $V_L$ overlapping in a face of smaller dimension. By induction, we are 
finished.

Next, suppose $2d+1\ge n$. We need to check that every valid regular subdivision can be refined to one 
whose cell has dimension $n$. If the two spanning facets are both simplices, then we can refine this 
regular subdivision to a regular triangulation, whose cell has $n$ degrees of freedom. If not, then as 
in the first case we can remove one of the overlap points from one facet, eventually reducing to the 
case where both facets are simplices. This completes the proof of the proposition.
\end{proof}

However, the second tropical secant variety (or complex) is not pure even for points in general 
position, and the 
first tropical secant variety is not pure if we do not assume that the points are in general position. 
We will give examples of these deviant cases, as well as an example of a non-generic linear subspace 
where the dimension is wrong, in Section~\ref{examples}.

\section{The one-dimensional case}\label{1-dim}
In this section, we apply the results of Section~\ref{maintheorem} to the one-dimensional case, i.e. 
when 
$L$ 
is a line in $\TP^n$. This corresponds to the case of a (projective) toric curve in $n$-space. In this 
case, we can completely compute the $k$-th tropical secant complex.

\begin{thm}
Let $L$ be a line in $\TP^n$, generated by $(r_1,\ldots,r_{n+1})$. Then the $k$-th secant complex
of $L$ 
consists of the set of lower faces of the cyclic polytope 
$C(2k,d-2)$ (i.e. $d-2$ points in dimension $2k$), where $d$ is the number of distinct elements of 
$\{r_1,\ldots,r_{n+1}\}$.
\end{thm}

\begin{proof}

The corresponding point configuration $V_L$ is $n+1$ points on the real line, located at
$r_1,\ldots,r_{n+1}$; we are looking for the space of height vectors $(x_0,\ldots,x_n)$ such that the
points $\{(r_1,x_1), \ldots, (r_{n+1},x_{n+1})\}$ have $k$ facets in their upper envelope which
together contain all the points. To begin with, if any $r_i$ and $r_j$ are identical, then the
corresponding $x_i$ and $x_j$ must also be, as otherwise whichever is lower will have that point not be
in the upper envelope at all. This reduces us to the case where the $r_i$ are distinct. We will also
assume that $r_0<r_1<\ldots<r_n$.

\begin{figure}
\begin{center}  \includegraphics{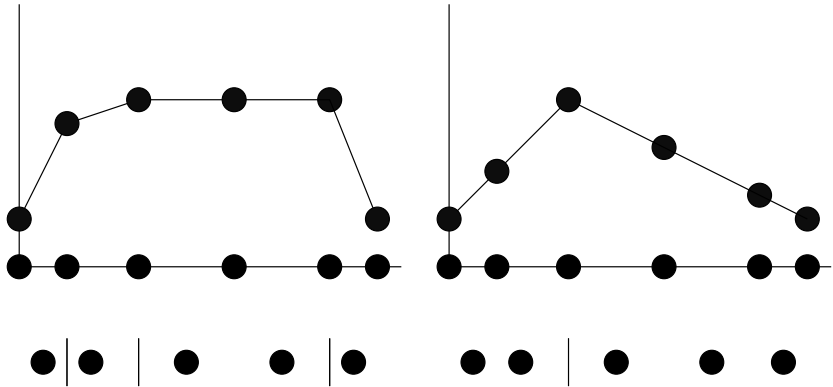}\end{center}
\caption[]{\label{1-cyclic}
Construction of bar patterns from height vectors of a one-dimensional point configuration.
}
\end{figure}

A height vector $(x_0,\ldots,x_n)$ is completely determined by the value
of $x_0$ and the slopes $s_i = \frac{x_{i+1}-x_i}{r_{i+1}-r_i}$ for $1\le
i\le n$. The condition that all of the points be in the upper envelope
reduces to the constraint $s_i\le s_{i-1}$ for $1\le i\le {n-1}$. When
$s_i = s_{i-1}$, the points $(r_{i-1}, x_{i-1})$, $(r_i,x_i)$, and
$(r_{i+1},x_{i+1})$ are all in the same facet of the corresponding regular
triangulation. We represent the $n$ slopes by points, putting a bar
between two points if the corresponding slopes differ; see
Figure~\ref{1-cyclic} for the corresponding illustration.

It is plain to see that the vertices of the complex provided by Corollary~\ref{complex} (for any $k$),
which are the coarsest nontrivial regular subdivisions, correspond to patterns with just one bar, of 
which there are $n-1$. 

Therefore, to compute the complex of the $k$-th secant variety, it suffices to compute the finest
regular subdivisions with $k+1$ facets whose union is all of the points. The $n+1$ lifted points are
represented by the $n+1$ spaces in between the dots (counting the beginning and the end), and the
facets of a lifted configuration are given by the closed intervals between adjacent bars. 

Suppose we have a bar pattern which is part of the $k$-th secant variety. The $k+1$ facets given by
Theorem~\ref{thm-main} clearly must contain the interval from the first 
space to the first bar. If the
space following the first bar is not filled by a bar, then in order for the union of these facets to
contain this space, we must also take the facet consisting of the interval from the first bar to the
second bar. If this is the case, then we can insert a bar in the space following the first bar, and
take the shortened facet from the new bar to the old second bar along with the other $k$ original
facets; we have just demonstrated that this augmented bar pattern still has the property.

Therefore, any bar pattern which corresponds to a facet of the tropical secant complex must have its 
first two bars adjacent. By identical reasoning, it in fact follows that the bars must come in adjacent 
pairs; the facets we take consist of every other interval, including the first and last. Counting 
carefully, we conclude that there must be $k$ pairs of consecutive bars or $2k$ bars in all.

So, the condition to be a facet of the tropical secant complex is that the included points must come in 
consecutive pairs. This is precisely the same as Gale's evenness condition~\cite{Z} for being a 
facet 
of the cyclic polytope $C(2k,n-1)$, except that our condition has the added stipulation that the 
initial and final segments of bars must have even length. This stipulation is easily seen to correspond 
exactly to the statement that the facet in question is a lower one, and so the $k$-th tropical 
secant complex is 
(combinatorially) isomorphic to the complex of lower faces of the cyclic polytope $C(2k,n-1)$ as 
desired.
\end{proof}

In the one-dimensional case, the $k$-th tropical secant complex is the same regardless of the (generic) 
linear subspace
chosen; this is because the oriented matroid of any $n$-point configuration in general position is the
same, and so the space of regular triangulations of these configurations are all combinatorially
isomorphic. In the two-dimensional case, corresponding to toric surfaces, the oriented matroid
statement is of course no longer true, and this leads to a wide variety of behavior of these tropical
secant varieties. We investigate this behavior in the next section.

\section{Examples}\label{examples}

In the previous section, we completely computed the tropical secant complexes of all one-dimensional
linear subspaces. In particular, these $k$-th secant complexes were all pure and contractible. In this
section, we compute a variety of examples, in which we show that these phenomena were purely
low-dimensional: even in two dimensions and even for $k=1$, the phylum of tropical secant complexes is
diverse, including species which are not pure and not contractible. We also demonstrate 
non-generic examples where the dimension of the tropical secant variety is not in accordance with 
Proposition~\ref{cxdim}. Rather than give the details of all of the computations, we present diagrams 
of the complexes associated to various configurations (with labeled vertices and facets) as well as pointing out some 
relevant aspects.

\begin{figure}
\begin{center}\includegraphics{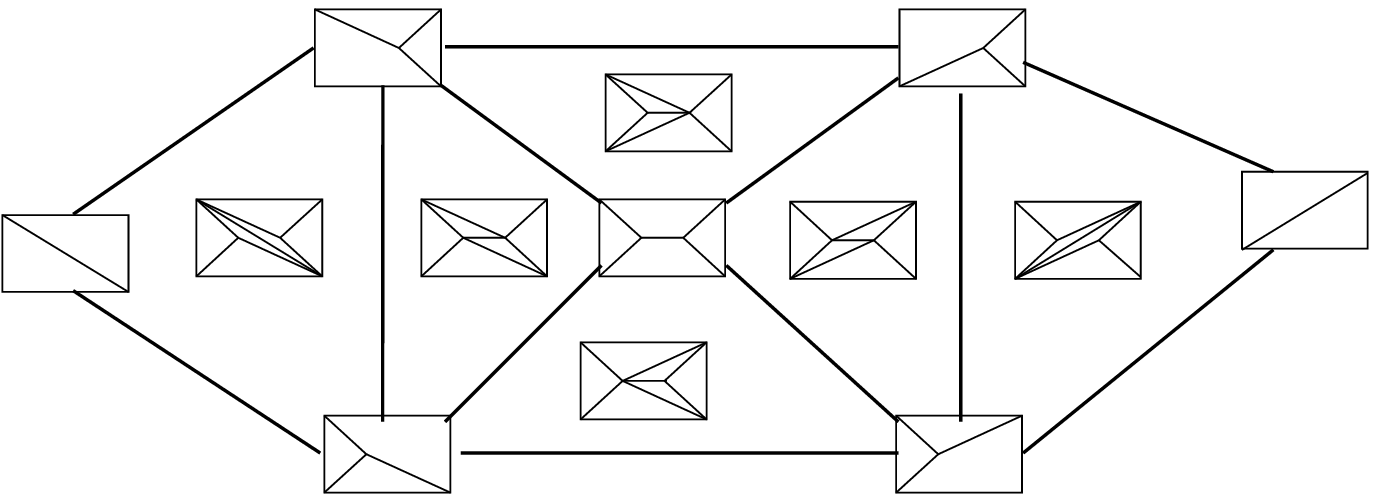}\end{center}
\caption[]{\label{square-plus-2}
First secant complex of $\left(\begin{array}{cccccc}
0 & 0 & 1 & 1 & 2 & 2 \\
0 & 3 & 1 & 2 & 0 & 3 
\end{array}
\right)$.}
\end{figure}

Our first example is a point configuration in general, non-convex position: a square with two points inside it. 
The corresponding linear subspace is generated by the rows of the matrix with the points as its columns:
\[
\left(\begin{array}{cccccc}
0 & 0 & 1 & 1 & 2 & 2 \\
0 & 3 & 1 & 2 & 0 & 3
\end{array}
\right).
\]
This linear subspace, a two-dimensional linear subspace in $\TP^5$, exhibits fairly normal behavior. Its first secant 
variety is the right dimension (five); indeed, the first secant variety is equal to the $k$-th secant 
variety for all 
$k>1$, since every regular subdivision of this point configuration has two facets whose union is the whole space 
(one containing the left three points and one containing the right three points.) This first tropical secant complex 
is shown in Figure~\ref{square-plus-2}. Since the points are in general position, it is pure by 
Proposition~\ref{first-pure}, and it is also contractible.

\begin{figure}
\begin{center}  \includegraphics{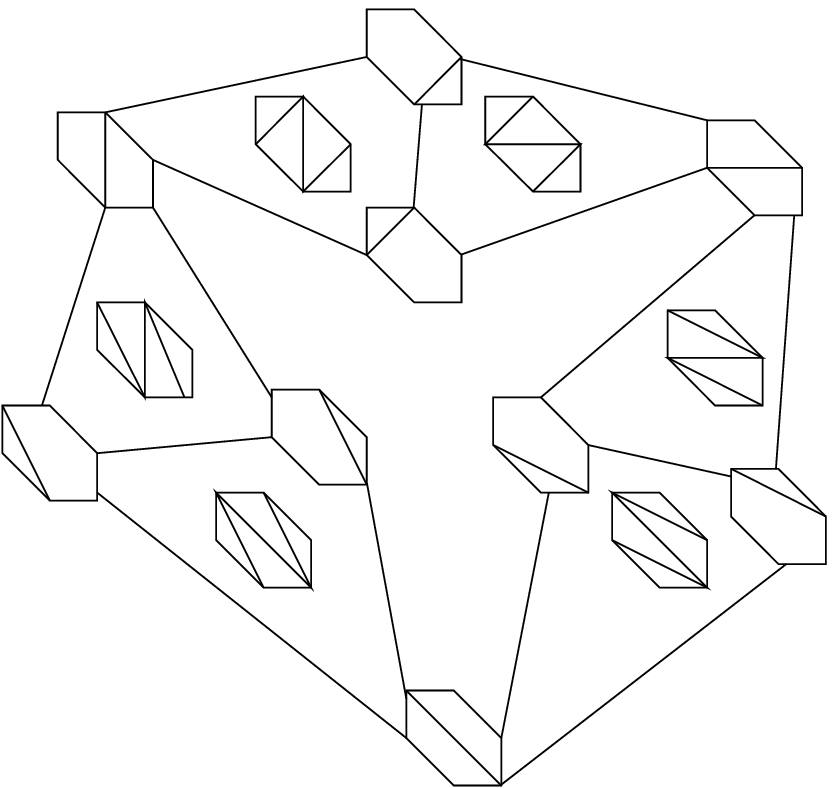}\end{center}
\caption[]{\label{n-gon}
First secant complex of 
$\left(\begin{array}{cccccc}
0 & 0 & 1 & 1 & 2 & 2 \\
0 & 1 & 0 & 2 & 1 & 2
\end{array}
\right)$.
}
\end{figure}

Our next example is $n$ points in general convex position, i.e. an $n$-gon. For $n=6$, the second secant complex is 
simply the dual of the secondary polytope of the hexagon, since every regular subdivision contains three facets whose 
union is all of the points. The first secant complex, however, is interesting: in particular, it is not contractible. 
For a general $n$-gon, the first secant complex's vertices will correspond to diagonals; the facets (it is a 
two-dimensional  complex) correspond to subdivisions given by three diagonals $\{(a,b+1), (a,b), (a+1, b)\}$. For 
$n\ge 7$, this complex forms a M\"{o}bius strip. For $n=6$, we show the diagram in Figure~\ref{n-gon}; this 
is the first secant complex of the linear subspace generated by the rows of 
\[
\left(\begin{array}{cccccc}
0 & 0 & 1 & 1 & 2 & 2 \\
0 & 1 & 0 & 2 & 1 & 2
\end{array}
\right).
\]
It is a subcomplex of the dual of the secondary polytope, which has 14 facets corresponding to the 14 triangulations 
of a hexagon. Six of these have the property that they have two facets whose union contains all of the vertices; 
these are the six facets in this first secant complex. 

\begin{figure}
\begin{center}  \includegraphics{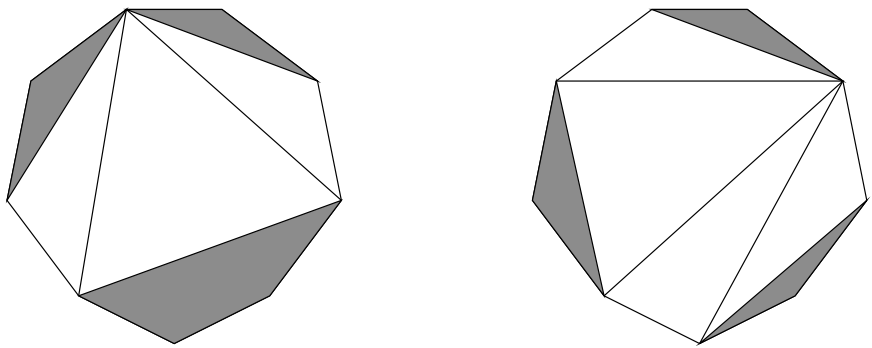}\end{center}
\caption{\label{impure-second}
Regular subdivisions corresponding to maximal faces of different dimension
in the second tropical secant complex of the linear subspace corresponding
to a nonagon.
}
\end{figure}

It is worth noting that the second secant complex of an $n$-gon is already not pure. Take $n=9$; then the regular 
subdivision on the left-hand side of Figure~\ref{impure-second} is a face of dimension four in the complex not 
contained in any face of dimension five, while the subdivision on the right-hand side is a face of dimension five. 
The spanning facets are shaded.

\begin{figure}
\begin{center}  \includegraphics{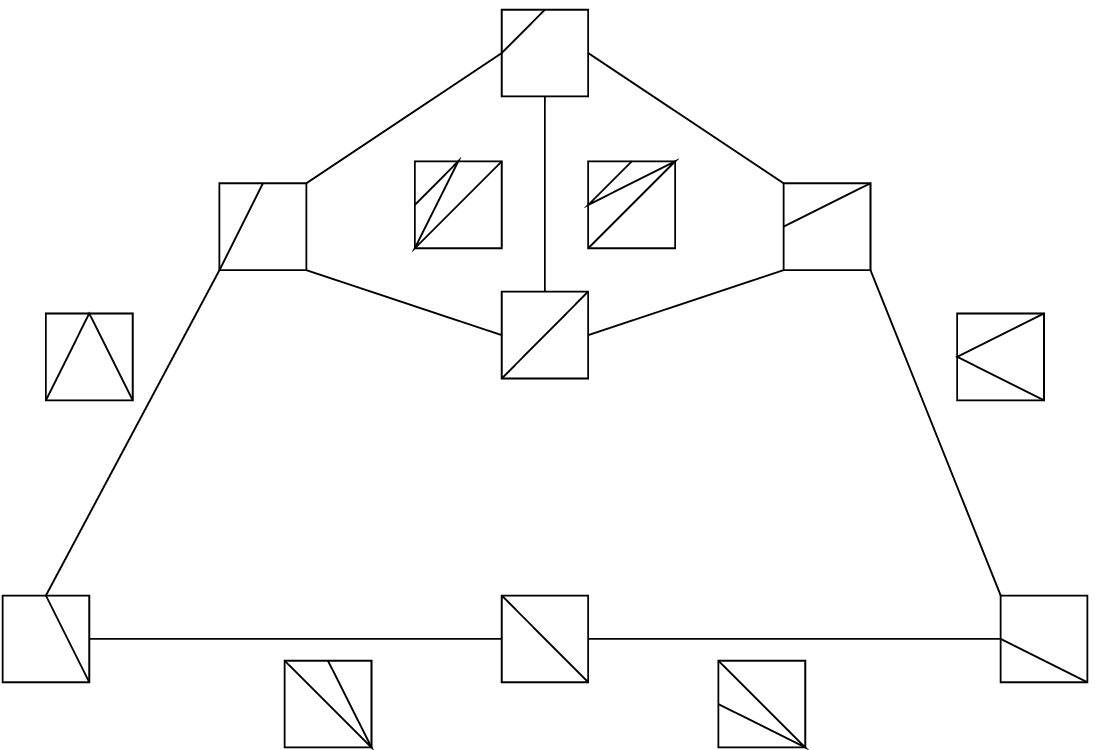}\end{center}
\caption[]{\label{impure-first}
First secant complex of $\left(\begin{array}{cccccc}
0 & 0 & 0 & 1 & 2 & 2 \\
0 & 1 & 2 & 0 & 0 & 2
\end{array}
\right)$.
}
\end{figure}

In our third example, the first secant complex is not even pure. The point configuration $V_L$ consists of two sets 
of three collinear points intersecting at a common end point, along with a point in the cone formed by the emanating 
rays, i.e. the point configuration corresponding to the linear subspace/matrix
\[
\left(\begin{array}{cccccc}
0 & 0 & 0 & 1 & 2 & 2 \\
0 & 1 & 2 & 0 & 0 & 2
\end{array}
\right).
\]

This complex consists of two triangles along with four edges, a "diamond ring" graph. It is
topologically equivalent to a circle and is decidedly impure. We can achieve the same behavior in
convex position in one higher dimension by taking a point configuration consisting of a cube with 
a point beyond a facet. The subdivision given by two opposite square pyramids including the extra
point and a
triangulation elsewhere is then a maximal face in the complex, but it has lower dimension than another
face of the complex, namely the one consisting of the extra point and three points of the facet nearest 
it, the square pyramid formed by the other five points, and simplices to fill out the remainder of the 
cube.

In both the second and third examples, the property responsible for the impureness is an oriented
matroid one. In the maximal secant complex face of inappropriate dimension, two facets of the spanning
set overlap, and ordinarily (if the points are in general position) we are able to deal with this by
removing a point from one of them, thus refining the subdivision. However, in the deficient cases, the
overlap point is an isthmus in the matroids of both spanning facets, and thus we cannot remove it from
either without degenerating them.

\begin{figure}
\begin{center}  \includegraphics{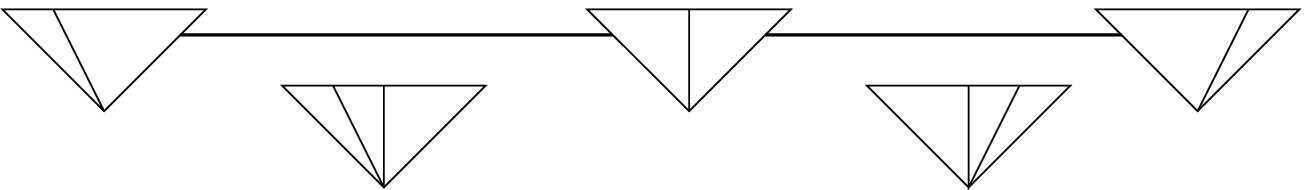}\end{center}
\caption[]{\label{five-on-a-line}
First secant complex of $\left(\begin{array}{cccccc}
0 & 0 & 0 & 0 & 0 & 2 \\
0 & 1 & 2 & 3 & 4 & 0
\end{array}
\right)$.
}
\end{figure}

The previous example shows that first secant complexes can be impure, although in those examples they at least have 
the right dimension. However, this need not be the case, as the next example shows. Here, the point configuration 
$V_L$ consists of one point off a line and five points on a line, corresponding to the linear subspace generated by
the rows of
\[
\left(\begin{array}{cccccc}
0 & 0 & 0 & 0 & 0 & 2 \\
0 & 1 & 2 & 3 & 4 & 0
\end{array}
\right).
\]
Proposition~\ref{cxdim} tells us that the first tropical secant variety should have dimension $(k+1)(d+1)-1 = 5$. 
However, it has only dimension four: that is, the first tropical secant complex is one-dimensional, not 
two-dimensional. Its vertices correspond to the coarsest possible nontrivial subdivisions, those with one dividing 
line segment; there are no two-faces, which correspond to having three dividing line segments, since the only regular 
subdivision with three dividing line segments does not have two facets whose union is the whole space. The first 
tropical secant complex, which consists of two line segments intersecting at a single point, is shown in 
Figure~\ref{five-on-a-line}.

This diversity among tropical secant varieties, even for first secant varieties of two-dimensional linear subspaces, 
is just a peek at the inviting theory. These complexes need not be pure, of the predicted dimension, or contractible. 
In the next section, we shift our focus to applying the theory to an already existing setup, namely the complex of 
matrices of fixed Barvinok rank.

\section{Barvinok rank}\label{barvinok}

The following definition (taken from~\cite{DSS}) is important in combinatorial optimization.

\begin{defn}
A matrix has {\em Barvinok rank} $k$ if it can be expressed as the tropical sum of $k$ matrices of 
tropical rank one, but not as the tropical sum of $k-1$ such matrices. A matrix $M$ has tropical rank 
one if we can write $M_{ij} = x_i + y_j$ for some $x_i$'s and $y_j$'s.
\end{defn}

Barvinok, Johnson, and Woeginger~\cite{BJW} showed that for matrices of fixed Barvinok rank $k$, the 
traveling salesman problem can be solved in polynomial time. Thus, an algorithm for finding the 
Barvinok rank of a matrix, or a description of the space of matrices of Barvinok rank $k$, is 
important. We have previously considered this problem in the papers~\cite{DSS} (with Bernd Sturmfels 
and Francisco Santos) and~\cite{DPts}.

Using the terminology in this paper, an attractive reformulation of the problem emerges. It is evident
from the definition that a matrix $M\in \RR^{d\times n}$ has Barvinok rank $k$ if it lies in the
$(k-1)$-st secant variety of $L$, where $L$ is the space consisting of all matrices of tropical rank
one. This is a linear subspace in $dn$ variables, defined by the equations
$M_{ij}+M_{kl}=M_{il}+M_{kj}$. A basis for this linear subspace of $\TP^{n\times d}$ is given by the
matrices $R_1,\ldots,R_{d-1}, C_1,\ldots,C_{n-1}$, where $R_i$ has 1's in the $i$-th row and 0's
everywhere else, and $C_j$ has 1's in the $j$-th column and 0's everywhere else. The columns of the
matrix with these rows are all 0-1 vectors with zero or one 1's among the first $d$ coordinates, and
zero or one 1's among the last $n-1$. This point configuration consists of the vertices of the product
of simplices $\Delta^{d-1}\times \Delta^{n-1}$.

Therefore, to determine the Barvinok rank of a matrix $M\in \RR^{d\times n}$, it suffices to consider
the regular triangulation induced by the corresponding height vector on $P=\Delta^{d-1}\times
\Delta^{n-1}$, i.e. by $\omega(i,j) = M_{ij}$, where $i\in [d]$ and $j\in [n]$. The Barvinok rank will 
be the smallest number of facets needed to cover all of the vertices of $P$. Similarly, the complexes 
of Barvinok rank $2,\ldots,\text{min}(d,n)$ are nested subcomplexes of the secondary polytope of $P$.

\end{document}